\documentclass[12pt]{amsart}
\usepackage[margin=1in]{geometry}
\usepackage{newtxtext}
\usepackage{newtxmath}
\usepackage[alphabetic]{amsrefs}
\usepackage{graphicx}
\usepackage{caption}
\usepackage{tikz}
	\usetikzlibrary{decorations.pathreplacing}
	\usetikzlibrary{decorations.pathmorphing}
	\usetikzlibrary{patterns}
\usepackage{setspace}
\usepackage[hidelinks]{hyperref}

\newcommand{\cd}{\cdot}
\newcommand{\ra}{\rightarrow}
\newcommand{\pr}{\prime}

\newcommand{\N}{\mathbb{N}}

\newcommand{\R}{\mathbb{R}}

\newtheorem{theorem}{Theorem}

\newtheorem{proposal}{Proposal}
\theoremstyle{definition}

\theoremstyle{remark}

\begin{document}

\title{The Proof of Kolmogorov-Arnold May Illuminate NN Learning}
\date{\today}

\author{Michael H. Freedman}
\address{\hspace{-\parindent}Michael H. Freedman, Center of Mathematical Sciences and Applications, Harvard University}
\email{freedmanm@google.com, mfreedman@cmsa.fas.harvard.edu}

\begin{abstract}
	Kolmogorov and Arnold, in answering Hilbert's 13th problem (in the context of continuous functions), laid the foundations for the modern theory of Neural Networks (NNs). Their proof divides the representation of a multivariate function into two steps: The first (non-linear) inter-layer map gives a universal embedding of the data manifold into a single hidden layer whose image is patterned in such a way that a subsequent dynamic can then be defined to solve for the second inter-layer map. I interpret this pattern as ``minor concentration'' of the almost everywhere defined Jacobians of the interlayer map. Minor concentration amounts to sparsity for higher exterior powers of the Jacobians. We present a conceptual argument for how such sparsity may set the stage for the emergence of successively higher order concepts in today's deep NNs and suggest two classes of experiments to test this hypothesis.
\end{abstract}

\maketitle

The Kolmogorov-Arnold theorem (KA) \cites{kol56,arn57} resolved Hilbert’s 13th problem in the context of continuous functions (Hilbert’s literal statement). It took some time for applied mathematicians and computer scientists to notice that it addresses the representation power of shallow, but highly non-linear, neural nets \cite{hed71}. The relevance of KA to machine learning (ML) has been much debated \cites{gp89,v91}.

The primary criticism is that the activation functions are not even first differentiable (even when the function $f$ to be represented is); indeed, they must be quite wild, and appear impossible to train. I agree with this criticism, but will argue here that it misses a more important point. The discussion of KA relevance has revolved around its statement. But in mathematics, proofs are generally more revealing of power than statements; I believe this particularly true in the case of KA. The purpose of this note is to call attention to some wisdom embedded within the proof which I believe will be useful for the training of neural nets (NNs). As far as resurrecting KA from the dustbin of ML history, this has already been done earlier this year in \cite{ZWV24}, where the original 1-hidden layer of KA has been deepened to many layers while the interlayer maps somewhat tamed (the authors use cubic $b$-splines).

Before explaining KA’s insight and its potential applicability, let me give a modern, optimized statement of KA quoted from \cite{mor21}. This statement is perhaps unnecessarily succinct from the NN perspective: it gets by with a single outer function $g$, not $2n+1$ outer functions which might appear more natural. Also, in the NN context one would expect the inner functions to carry a second index and not merely be rationally independent scalars times a single index function $\phi_k$. But we take this statement as representative and refer the reader to Morris for more historical information on the individual contributions.

\begin{theorem}[Kolmogorov, Arnold, Kahane, Lorentz, and Sprecher]
	For any $n \in \N$, $n \geq 2$, there exist real numbers $\lambda_1, \lambda_2, \dots, \lambda_n$ and continuous functions $\phi_k: \mathbb{I} \ra \R$, for $k = 1, \dots, m$, $m \geq 2n+1$, with the property that for every continuous function $f: \mathbb{I}^n \ra \R$ there exists a continuous function $g: \R \ra \R$ such that for each $(x_1,x_2,\dots,x_n) \in \mathbb{I}^n$,
	\[
	f(x_1,\dots,x_n) = \sum_{k=1}^m g(\lambda_1 \phi_k(x_1) + \cdots + \lambda_n \phi_k(x_n)).
	\]
\end{theorem}

\begin{figure}[ht]
	\centering
	\begin{tikzpicture}[scale=1.1]
		\draw[fill=black] (-4,1) circle (0.2ex);
		\node at (-4.4,1) {$x_1$};
		\draw[fill=black] (-4,-1) circle (0.2ex);
		\node at (-4.4,-1) {$x_2$};

		\draw (-4,1) -- (0,2);
		\draw (-4,1) -- (0,1);
		\draw (-4,1) -- (0,0);
		\draw (-4,1) -- (0,-1);
		\draw (-4,1) -- (0,-2);

		\node at (-2,1.8) {$\phi_{p,q}$};
		
		\draw[dashed] (-4,-1) -- (0,0);
		\draw[dashed] (-4,-1) -- (0,-1);
		\draw[dashed] (-4,-1) -- (0,-2);
		\draw[dashed] (-4,-1) -- (0,1);
		\draw[dashed] (-4,-1) -- (0,2);

		\draw[fill=black] (0,2) circle (0.2ex);
		\draw[fill=black] (0,1) circle (0.2ex);
		\draw[fill=black] (0,0) circle (0.2ex);
		\draw[fill=black] (0,-1) circle (0.2ex);
		\draw[fill=black] (0,-2) circle (0.2ex);

		\node at (2,1.3) {$g$};

		\draw (0,2) -- (4,0);
		\draw (0,1) -- (4,0);
		\draw (0,0) -- (4,0);
		\draw (0,-1) -- (4,0);
		\draw (0,-2) -- (4,0);

		\draw[fill=black] (4,0) circle (0.2ex);
		\node at (4.9,0) {$f(x_1,x_2)$};
	\end{tikzpicture}
	\caption{}
\end{figure}

The proof of KA is divided into the construction of inner functions $\phi_i$ and an outer function $g$. The inner functions taken together constitute an embedding, which we now denote $\Phi: \mathbb{I}^n \ra \R^m$, $m \geq 2n+1$. We will append a subscript $\Phi_i$ when referring to a neuron coordinate of the embedding. $\Phi$ is universal; it can be chosen once and for all, independent of the function $f$ to be represented. Our chief lesson regards $\Phi$, although I also will make a comment on the non-linear dynamic used to converge to $g$ (given $f$).

In the mathematically idealized case of KA, feasible $\Phi$ are actually dense in the space of continuous maps $\mathcal{C}(\mathbb{I}^n,\R^m)$; they also have a crucial local feature, about to be described. This fact is reminiscent of the observation (for example see \cites{gmrm23,skmf24}) that when NNs are randomly initialized, training is often seen to merely perturb the values of the early layers. These layers seem to be more in the business of assuming a form conducive to the learning of later layers, rather than attending to the particular data themselves.

The local property crucial to $\Phi$ is a kind of irregular staircase structure as might be used to approximate the general continuous function by one which is piecewise constant. More specifically, $\Phi$ may be taken to be Lipschitz \cite{ak17}, and so by Rademacher’s theorem will be differentiable a.e. Such $\Phi$ constructed in the proof of KA will have the property that at every point $x \in \mathbb{I}^n$ where the Jacobian is defined it will, as an $n \times m$ matrix in the ``neuron basis,'' have $m - n$ columns consisting entirely of zeros. So only one of its ($m$ choose $n$) $n \times n$ minors can be nonzero. (Of course, which minor is active will vary for with $x$.) To picture what this means, imagine an irregular stack of sugar cubes of different sizes but with all faces perpendicular to one of the three coordinate axes, $x$, $y$, or $z$. Another condition says that the ``inactive'' coordinate on the $w$-perpendicular faces ($w = x$, $y$, or $z$) must take distinct values. The actual situation is more like an irregular structure of $n$-D cubes wafting through $m$-D space, now with $m-n$ inactive directions at any generic point, and all these inactive coordinate values distinct. Being distinct is what sets the stage towards finding at least an approximate outer function $g_0$. The $m-n$ inactive directions are crucial to the construction of $g_0$, and it is important that $(m-n)>n$, so that the inactive coordinates are in the majority at all generic points of $\Phi(x)$.

The very rough idea for how to choose $g_0$ is to use the $(m-n)$ inactive function value on each ``sugar cube face'' as a hint: If that value $f(x)$ is uniformly positive (negative) on the face with inactive coordinate $\Phi_i(x)$, give $g_0(\Phi_i(x))$ a small constant positive (negative) value and extend in a convex-PL manner. In the regime $(m-n)>n$, it may be checked that such a strategy produces a $g_0$ leading to a useful approximation:
\[
\Big\lVert f(x_1,\dots,x_n) - \sum_{i=1}^n g_0(\Phi_i(x_1,\dots,x_n)) \Big\rVert < \mathrm{const} \lVert f(x_1,\dots,x_n) \rVert,\ \mathrm{const} < 1.
\]

Now iterating on the approximation error (in a manner reminiscence of Resnets), and relying on additivity at the output node, one produces a series $g_0, g_1, \dots$ which converges exponentially fast to the desired outer function $g$ for which the approximation error has been driven to zero. Many inactive directions at each point are essential since they provide a stationary coordinate value on which a guess for $g_0$ can be made: A small positive constant if $f(x_1, \dots ,x_n)$ is positive throughout the face and a small negative constant if $f$ is negative throughout the face.

This sketch of the construction of the outer function g reveals the importance of the stationarity condition and its implication of vanishing Jacobian minors in the embedding $\Phi$. Remembering that admissible $\Phi$ are dense, we see that the job of the first layer is simply to impress a certain microscopic texture on $\Phi$ which makes g constructable (if not learnable). This is the division of labor that feels so striking. $\Phi$ does not try to learn anything but sets up for success the task of finding $g$.

I wonder if a similar division of labor occurs between early and later stages of a NN as it is trained, and propose:

\begin{proposal}\label{prop1}
	Search the natural \emph{Jacobian maps} for \emph{minor concentration}.
\end{proposal}

First, what are these natural Jacobians, and second, what is minor concentration? One may consider either a conventional feedforward DNN or the recently proposed KANets \cite{ZWV24}, or more elaborate architectures with residual connections and self-attention. The idea is to watch the data manifold flow through the NN. Initially at $t=0$, the data manifold is $\mathbb{I}^n$, the values stored in the $n$ input neurons, and at time $j$ when the input has been imaged in the $j$th layer (say after application of the activation functions). Let $k_{ij}$, $0 \leq i,j < \ell_{\text{max}}$, be the map between layers $i$ and $j$, with $i=0$ being a very interesting special case. At any point $y$ in the evolving data manifold ($y=x$, the previously defined coordinate on $\mathbb{I}^n$, if $i=0$) we may consider $J_y(k_{ij}) := J_{ij}y$, the Jacobian matrix (in neuron coordinates) of the forward map $k_{ij}$ at point $y$. The proposal is to find some context where it is possible to search (or at least sample) this huge family of Jacobians and look for minor concentration. Based on understanding KA one may expect, particularly in some initial segment of the NN (perhaps all the way up to the penultimate layer), that training will cause these Jacobians to have their minors concentrated far beyond what would be seen in a random collection (e.g.\ Gaussian distributed entries) of linear maps. It seems best at first to start from the beginning and set $i=0$. That is, one should look ``in the wild'' to see if a KA-style of data preparation is discovered, at least for certain tasks, and certain architectures, as a result of training.

\subsection*{What is minor concentration?}
Given a $p \times q$ matrix $M$ and $h$, $0< h \leq \min\{p,q\}$, minor concentration should be some quantity that measures how far from uniform the distribution of the absolute value of the $\binom{p}{h} \binom{q}{h}$ $h \times h$ minors of $M$ are. The easiest formula is:
\[
\operatorname{MC}(M) = \frac{\lVert \text{absolute value of minors} \rVert_2}{\lVert \text{absolute value of minors} \rVert_1}.
\]
That is, the ratio of the $L_2$ and $L_1$ norms.

As will become clear in the thought experiment below, one might also wish to design special purpose minor concentration functionals which reward not just a large, isolated minor but respond to a row or column of large minors inside $\Lambda^h M$, where the superscript $h$ indicates the $h$th exterior power. This type of pattern for large and small minors could plausibly arise through training.

\begin{proposal}
	Evaluate the effect of \emph{forced minor concentration} on \emph{learning}.
\end{proposal}

Independent of the results of Proposal \ref{prop1}, one could study the effect of encouraging/discouraging minor concentration during training. One way to do this would be to alternate traditional training protocols with a novel step where a layer map would be updated along the gradient of a new objective function, one measuring MC. If MC is defined as above, these interleaved steps could be $\delta^\pr \cd \nabla \mathrm{MC}_y$. That is, with MC-learning rate $\delta^\pr$, move a layer map in the direction of increased minor concentration. The gradient might be taken over the last layer-map variables---although there are many choices with which to experiment. One might expect to find some portions of the NN and some stages of training where ``concepts are gelling'' (see below) in which a positive $\delta^\pr$ would accelerate learning and a negative $\delta^\pr$ delay learning. Interestingly, there could also be regimes where the reverse applied: At an early stage of learning, it could be harmful to have a concept gel prematurely, the NN may be better off keeping an open mind (Zen-like) mind. Experiments hindering/enhancing MC could get to the heart of ``how learning works'' \cite{gmrm23}.

\subsection*{Thought experiment/example}
In the KA proof for $n=2$, $m = 2n+1 = 5$, the differential $J_x$ at a generic point $x$ of $\Phi$ will be a $2 \times 5$ matrix $M$, and the salient feature of $J_x$ is that for any given $x \in \mathbb{I}^2$, $M$ will have three of its five columns entirely zero (which three are zero will vary as $x$ varies), so in particular the $\Lambda^2 M$, the second exterior power, will be a $1 \times 10$ matrix with at most a single nonzero entry. This is the prototypical example of minor concentration. Let’s consider another.

Say that between layers $j$ and $j+1$, the Jacobian of the non-linear map is $J_y = M$. In this thought experiment, let’s imagine we have a convolutional NN analyzing a picture of a face. Suppose the $j$th layer has learned some lines encoded in four neurons, determining a vector in $\R^4$, and say the $(j+1)$th layer has three neurons spanning in $\R^3$ and may be looking to define a feature. Now $M$ is a $4 \times 3$ matrix. Arbitrarily pick $h=2$, meaning we will be inquiring on how two ``line'' neurons stimulate two ``feature neurons.'' Now $\Lambda^2 M$ is a $6 \times 3$ matrix. Suppose we see just one of its three rows with large values while the other 12 values are close to zero. Here is a story that we might tell: Let the feature-neurons $a$, $b$, and $c$ label the rows of $M$, so $ab$, $bc$, and $ca$ label the rows of $\Lambda^2 M$. Suppose it is the $bc$ row with large minors. This means as we consider the six pairs of ``line neurons,'' as their information varies to first order, we see the information in the $bc$ pair vary robustly, but not so with the $ab$ and $ca$ pairs. That would be consistent with each of these three pairs being responsible for detecting a feature, say, eye, nose, and mouth respectively. The most parsimonious explanation for the muted response of the eye and mouth pairs would be that the lines are not suggesting either of these features, so varying the incoming line data leaves the $ab$ and $ca$ neurons in their ``not an eye'' and ``not a mouth'' base point states respectively. There is little or no variation in the eye or mouth recognizing neuron pairs because they are not finding what they have been trained to look for. On the other hand, at least some of the six pairwise line characteristics seem to be describing a variety of nose shapes as they move about. Because training involves both forward- and back-propagation, similar scenarios can be constructed where a column of $\Lambda^h M$ would be heavy. Thus, in addition to the most na\"{i}ve measure of minor concentration $\operatorname{MC}(M)$ (above), one should also look for row-wise and column-wise concentrations in $\Lambda^h M$. Looking for such concentrations would always come down to comparing the appropriate functional on Jacobians from trained NNs with random matrices.

It is worth pointing out that when $h=1$, minor concentration is quite close to the notion of ``sparsity''; 1-minors concentrate when a few matrix entries are large (in absolute value) and the rest small. So, minor concentration for general $h$ is close to looking for sparsity in the exterior powers $\Lambda^h M$. It is worth commenting that the old notion of sparsity and its generalization to k-minor concentration only makes sense in the context of linear maps between based vector spaces. Without preferred bases, neither the individual entries in matrices, nor their exterior powers, have any invariant meaning.

\subsection*{A caveat}
I’d like to mention a final point which arose in conversation with Boris Hanin relating to the phenomenon of neural collapse. Observations presented in \cite{koth22} suggest, for example, if a convolutional network is trained to recognize farm animals, say: cow, sheep, pig, horse, that the penultimate layer may have the data compressed to some regular tetrahedron, with the four animal types at the four vertices, but these tetrahedron vertices will not generally be oriented along neuron axes. To put them in this form a rotation may be required. This suggests information midway though the net may often be held diffusely and be difficult to detect within a small subset of neurons, and hence difficult to detect through minor concentration, at least for small $h$. However, sparsification, if applied during training, may already be emphasizing the neuron basis, and so facilitate concentration of $h$-minors for $h>1$. If this is so, it could present an interesting avenue towards understanding how sparsification at the original matrix level $M$ could enhance minor concentration and facilitate learning.

\bibliography{references}

\begin{bibdiv}
\begin{biblist}

\bib{ak17}{article}{
      author={Actor, Jonas},
      author={Knepley, Matthew~G.},
       title={An algorithm for computing {L}ipschitz inner functions in
  {K}olmogorov's superposition theorem},
        date={2014},
      eprint={arXiv:1712.08286},
}

\bib{arn57}{article}{
      author={Arnold, Vladimir},
       title={On functions of three variables},
        date={1957},
     journal={Proceedings of the USSR Academy of Sciences},
      volume={114},
       pages={679\ndash 681},
        note={English translation: American Mathematical Society Translations,
  ``Sixteen Papers on Analysis'' (1963), 51--54.},
}

\bib{gmrm23}{article}{
      author={Guth, Florentin},
      author={Ménard, Brice},
      author={Rochette, Gaspar},
      author={Mallat, Stéphane},
       title={A rainbow in deep network black boxes},
        date={2023},
      eprint={arXiv:2305.18512},
}

\bib{gp89}{article}{
      author={Girosi, Federico},
      author={Poggio, Tomaso},
       title={Representation properties of networks: {K}olmogorov's theorem is
  irrelevant},
        date={1989},
     journal={Neural Computation},
      volume={1},
      number={4},
       pages={465\ndash 469},
}

\bib{hed71}{inproceedings}{
      author={Hedberg, Torbj\"{o}rn},
       title={The {K}olmogorov superposition theorem},
        date={1971},
   booktitle={{Topics in Approximation Theory}},
      series={Lecture Notes in Mathematics},
      volume={187},
       pages={267\ndash 275},
}

\bib{v91}{article}{
      author={Kůrková, Věra},
       title={Kolmogorov's theorem is relevant},
        date={1991},
     journal={Neural Computation},
      volume={3},
      number={4},
       pages={617\ndash 622},
}

\bib{kol56}{article}{
      author={Kolmogorov, Andrey},
       title={On the representation of continuous functions of several
  variables by superpositions of continuous functions of a smaller number of
  variables},
        date={1956},
     journal={Proceedings of the USSR Academy of Sciences},
      volume={108},
       pages={179\ndash 182},
        note={English translation: American Mathematical Society Translations,
  ``Twelve Papers on Algebra and Real Functions'' (1961), 369--373.},
}

\bib{koth22}{article}{
      author={Kothapalli, Vignesh},
       title={Neural collapse: A review on modelling principles and
  generalization},
        date={2022},
      eprint={arXiv:2206.04041},
}

\bib{ZWV24}{article}{
      author={Liu, Ziming},
      author={Wang, Yixuan},
      author={Vaidya, Sachin},
      author={Ruehle, Fabian},
      author={Halverson, James},
      author={Soljačić, Marin},
      author={Hou, Thomas~Y.},
      author={Tegmark, Max},
       title={{KAN: Kolmogorov-Arnold} networks},
        date={2024},
      eprint={arXiv:2404.19756},
}

\bib{mor21}{article}{
      author={Morris, Sidney~A.},
       title={Hilbert 13: Are there any genuine continuous multivariate
  real-valued functions?},
        date={2021},
     journal={Bulletin of the American Mathematical Society},
      volume={58},
       pages={107\ndash 118},
}

\bib{skmf24}{article}{
      author={Shirai, Norihito},
      author={Kubo, Kenji},
      author={Mitarai, Kosuke},
      author={Fujii, Keisuke},
       title={Quantum tangent kernel},
        date={2024},
     journal={Physical Review Research},
      volume={6},
      number={3},
       pages={033179},
}

\end{biblist}
\end{bibdiv}

\end{document}